\documentclass{elsarticle}
 \pdfoutput=1
  \usepackage{amssymb,amsmath,latexsym,theorem,eepic,enumerate,
  graphicx,QED,shapepar}
  \usepackage[bookmarks=true,
              bookmarksnumbered=true, breaklinks=true,
              pdfstartview=FitH, hyperfigures=false,
              plainpages=false, naturalnames=true,
              colorlinks=true,pagebackref=true,
              pdfpagelabels]{hyperref}
  \usepackage{xypic,times}
  \xyoption{curve}
  \usepackage{geometry,color}
{\theorembodyfont{\rmfamily}\newtheorem{example}{Example}[section]}
{\theorembodyfont{\rmfamily}\newtheorem{Def}[example]{Definition}}
{\theorembodyfont{\rmfamily}\newtheorem{crit}[example]{Criterion}}
{\theorembodyfont{\rmfamily}\newtheorem{anomaly}[example]{Anomaly}}
{\theorembodyfont{\rmfamily}\newtheorem{strategic}[example]{Strategic Analysis}}
{\theorembodyfont{\rmfamily}\newtheorem{remark}[example]{Remark}}

\newtheorem{theorem}[example]{Theorem}
{\theorembodyfont{\rmfamily}}
\newtheorem{cor}[example]{Corollary}

  \def\FTop{\mathsf{FTop}}
 
  \def\Cub{\mathsf{Cub}}
  \def\Crs{\mathsf{Crs}}

  \def\Crs{\mathsf{Crs}}
  \def\XMod{\mathsf{XMod}}
  \def\Cub{\mathsf{Cub}}

\def\Crs{\mathsf{Crs}}

\def\Gpd{\mathsf{Gpd}}

\newcommand{\midsq}[1]{\save\go[0,0];[1,1]:(0.5,0)\drop{#1}\restore}

\def\midsq{\ar @{} [dr] |}

\newcommand{\threeaxes}[3]{\def\objectstyle{\scriptstyle}  \objectmargin={0pt}
\xy
(0,0)*+{}="a",(0,-6)*+{\rule{0em}{1.5ex}#2}="b",(7,0)*+{\;#1}="c",
(14,-3)*+{\;#3}="d" \ar@{->} "a";"b" \ar @{->}"a";"c"  \ar
@{->}"a";"d"\endxy }

\newcommand{\directs}[2]{\def\objectstyle{\scriptstyle}  \objectmargin={0pt}
\xy
(0,4)*+{}="a",(0,-2)*+{\rule{0em}{1.5ex}#2}="b",(7,4)*+{\;#1}="c"
\ar@{->} "a";"b" \ar @{->}"a";"c" \endxy }

\newcommand{\xdirects}[2]{\def\objectstyle{\scriptstyle}  \objectmargin={0pt}
\xy
(0,0)*+{}="a",(0,-6)*+{\rule{0em}{1.5ex}#2}="b",(7,0)*+{\;#1}="c"
\ar@{->} "a";"b" \ar @{->}"a";"c" \endxy }
\newcommand{\sdirects}[2]{\def\objectstyle{\scriptstyle}  \objectmargin={0pt}
\xy
(0,2.2)*+{}="a",(0,-2.5)*+{\rule{0em}{1.5ex}#2}="b",(7,2.2)*+{\;#1}="c"
\ar@{->} "a";"b" \ar @{->}"a";"c" \endxy }

\newcommand{\bl}{\mbox{\rule{0.08em}{1.7ex}\hspace{-0.00em}\rule{0.7em}{0.2ex}}}

\newcommand{\br}{\mbox{\rule{0.7em}{0.2ex}\hspace{-0.04em}\rule{0.08em}{1.7ex}}}

\newcommand{\tr}{\mbox{\rule[1.5ex]{0.7em}{0.2ex}\hspace{-0.03em}\rule{0.08em}{1.7ex}}}

\newcommand{\tl}{\mbox{\rule{0.08em}{1.7ex}\rule[1.54ex]{0.7em}{0.2ex}}}

\newcommand{\hh}{\mbox{\rule{0.7em}{0.2ex}\hspace{-0.7em}\rule[1.5ex]{0.70em}{0.2ex}}}

\newcommand{\vv}{\mbox{\rule{0.08em}{1.7ex}\hspace{0.6em}\rule{0.08em}{1.7ex}}}

\newcommand{\sq}{\mbox{\rule{0.08em}{1.7ex}\hspace{-0.00em}\rule{0.7em}{0.2ex}\hspace{-0.7em}\rule[1.54ex]{0.7em}{0.2ex}\hspace{-0.03em}\rule{0.08em}{1.7ex}}}

\newcommand{\tsq}{\mbox{\rule{0.04em}{1.55ex}\hspace{-0.00em}\rule{0.7em}{0.1ex}\hspace{-0.7em}\rule[1.5ex]{0.7em}{0.1ex}\hspace{-0.03em}\rule{0.04em}{1.55ex}}}

\def\rho{\varrho}

\def\B{\beta}

\def\epsilon{\varepsilon}

\def\B{\beta}

\def\epsilon{\varepsilon}

\def\bu{\bullet}
\def\colim{\mathrm{colim }}

 \def\colim{\mathop{\rm colim}\nolimits}

 \def\AD{\mathcal{AD}}

\def\ogpd{$\omega$-$\mathsf{Gpd}$}

\def\epsilon{\varepsilon}
\def\H{\mathbb H}
\def\B{\mathbb B}
\def\colim{\operatorname{colim}}

\def\xybiglabels{\def\labelstyle{\textstyle}}
\begin{document}
\begin{frontmatter}
 \title{Modelling and computing homotopy types: I \tnoteref{mytitlenote}}
 \tnotetext[mytitlenote]{To aqppear in a special issue of Idagationes Math. in honor of L.E.J. Brouwer to appear in 2017. This article, has developed from presentations given at Liverpool, Durham,  Kansas SU, Chicago, IHP, Galway, Aveiro, whom I thank for their hospitality. }

  \author{Ronald Brown}
  \address{School of Computer Science, Bangor University}

  \begin{abstract}
The aim of this article is to explain a philosophy for applying the 1-dimensional  and higher dimensional Seifert-van Kampen Theorems, and how the use of groupoids and strict higher groupoids resolves some foundational anomalies in  algebraic topology at the border between homology and homotopy. We   explain some applications to filtered spaces, and special cases of them, while a  sequel  will show the relevance to $n$-cubes of pointed spaces.
\end{abstract}

\end{frontmatter}


\tableofcontents

\section*{Introduction}\label{sec:intro}
\addcontentsline{toc}{section}{Introduction}
The usual algebraic topology approach to homotopy types is to aim first to calculate specific invariants and then look for some way of amalgamating these to determine a homotopy type. A standard complication is that in a space, identifications in low dimensions affect high dimensional homotopy invariants.

The philosophy which we state in outline in Section \ref{sec:phil} was developed by the author and others since about 1965 and was in response to a gradual realisation of certain anomalies in traditional algebraic topology which we explain in Section \ref{sec:anomal}.

A basic reason for formulating and trying to analyse these issues is given by Einstein in the correspondence published in \cite{Ein90}:

\begin{quote}
$\ldots$  the following questions must burningly interest me as a disciple of science: What goal will be reached by the science to which I am dedicating myself? What is essential and what is based only on the accidents of development? $\ldots$ Concepts which have proved useful for ordering things easily assume so great an authority over us, that we forget their terrestrial origin and accept them as unalterable facts.  $\ldots$  It is therefore not just an idle game to exercise our ability to analyse familiar concepts, and to demonstrate the conditions on which their justification and usefulness depend, and the way in which these developed, little by little $\ldots$
\end{quote}

One conclusion from progress made with resolving some of these anomalies is the difficulty of dealing with ``bare" topological spaces, i.e. spaces with no additional structure. In order to calculate something about a space $X$ we need some data on $X$ and this data will usually have some structure. Thus our algebraic invariants will need to recognise such structure. Further, in view of the influence of low dimensions on high dimensions, we need algebraic invariants which have structure in a range of dimensions. The  single base point which is necessary to define the $n$th homotopy group is a rather limited kind of structure.

A traditional view\footnote{This was expressed  by Saunders Mac Lane to Philip Higgins and the author in London in 1972.  } on a prospective  Seifert-van Kampen Theorem in dimension 2 was that we need to determine the Postnikov invariant $k_X \in  H^3(\pi_1(X),\pi_2(X))$ when $X$ is a union. So we need $\pi_1(X)$ by the standard Seifert-van Kampen Theorem; then find $\pi_2(X)$ of a union; then find $k_X$  in terms of the invariants of the pieces. The last two steps seem impossible.

However if $X$  is a CW-complex, then by the main result of \cite{MacW50}, $k_X$ is determined by the crossed module (see  Definition \ref{defxmod})
$$ \pi_2(X, X^ 1,x) \to  \pi_1(X^1, x).$$
 The advantage of this approach is that crossed modules form an algebraic structure for which colimits exist and can in principle be calculated; further they can be thought of as structures in dimensions 1 and 2.   So we can at least formulate the possibility that this crossed module of a union could be the colimit of the crossed modules of the pieces. What is much more difficult is to see how to prove, or disprove,  such a result.

A key to the Mac Lane-Whitehead paper was the result in \cite[Section 16]{W49CHII}, which we call Whitehead's Theorem,   that the crossed module
$$\pi_2(X \cup \{e^2_\lambda \},X, x) \to \pi_1(X, x)$$
could be described as a ``free crossed module''\footnote{This result is sometimes stated but rarely proved in texts on algebraic topology. Whitehead's proof involves ideas of transversality and knot theory. An exposition is given in \cite{B80}. } on the characteristic maps of the 2-cells $e^2_\lambda, \lambda \in \Lambda$ (see Diagram \eqref{eqn:freexmod}). Thus we see a universal property in $2$-dimensional relative homotopy theory. So the authors of \cite{BH78sec} began in 1974 to look for proving a 2-dimensional van Kampen theorem in a {\it relative}  situation, i.e. one with more structure; see the Strategic Analysis \ref{thm:strategy}. This led to the proof and formulation of a pushout  theorem, see Diagram  \eqref{equ:2push},  for Whitehead's crossed modules, \cite{BH78sec}.

The intuition that there might be a $2$-dimensional Seifert--van Kampen Theorem came in 1965 with an idea for the use of forms of double groupoids, although an appropriate generalisation of the fundamental groupoid was lacking.  We explain more on this idea in Sections \ref{sec:hdim}ff.

A further feature that seemed to me  of interest was  that the Seifert-van Kampen theorem is a {\it local-to-global theorem;} such theorems are of wide importance in mathematics and science. Further, the standard Seifert-van Kampen theorem for groups is a nonabelian theorem of this type, which is unusual. Algebraic models which could allow a higher dimensional version have the possibility of being really new. Such a view seemed therefore well worth pursuing, although it has been termed ``idiosyncratic". It can now be seen as related to issues\footnote{local to global, increase in dimensions, commutative to noncommutative, linear to non linear, geometry versus algebra} discussed in  \cite{Atiyah}. Any dichotomy between geometry and algebra as suggested in that paper is softened by the methods given in Section \ref{sec:nbdata}. Tools of category theory, and of combinatorial group theory,  are also needed for the effective application of the theory.

\section{	Statement of the philosophy}\label{sec:phil}
The general philosophy which has developed is to consider categories and functors
$$  \xybiglabels \xymatrix{ \left(\txt{Topological\\ Data }\right) \ar@<1ex> [r]^{\H}  &\ar @<1ex> [l] ^{\B} \left(\txt{Algebraic \\ Data }\right)}$$
such that:
\begin{crit}

$\H$ is homotopically defined.  \label{crit1}
\end{crit}

\begin{crit}
$\H \B$ is naturally equivalent to $1$. \label{crit2}
\end{crit}
\begin{crit}
The Topological Data has a notion of {\it connected}.  \label{crit3}
\end{crit}
\begin{crit}
For all Algebraic Data $\mathsf A$, the Topological Data $\B \mathsf A$ is connected.   \label{crit4}
\end{crit}

\begin{crit}
If the Topological Data $\mathsf X$ has underlying space $X$ which is the union of a family $\mathcal U= \{U_\lambda, \lambda \in \Lambda\}$ of open sets, and each corresponding Topological Data
$ \mathsf U_\nu$ is connected for all finite intersections $U_\nu$ of the sets of $\mathcal U$, then \\(a) $\mathsf X$ is connected, and (b) the following morphism induced by inclusions $U_\nu \to X$
\begin{equation*}
  \colim_\nu \H \mathsf U_\nu \to \H \mathsf X
\end{equation*}
is an isomorphism.
  \label{crit5}
\end{crit}

We will give detail on the notion of ``connected" in specific examples later (see Definitions \ref{def:connpair}, \ref{Def:connfiltered}): this condition does of course give a limitation on the applicability of such theorems, limitations which are familiar in algebraic topology. In fact considerable effort by many writers has gone into developing connectivity results, such as the theorems of Blakers-Massey, and of $n$-ad connectivity, sometimes proved by subtle homotopy theoretic arguments, or general position and transversality, see for example \cite{cubhom}. In practice, though, we would like an algebraic calculation, and these usually imply connectivity results, since often a kind of tensor product is involved, and this is trivial if one of the factors is trivial.
The aim is precise algebraic colimit calculations of some homotopy types.

\section{Broad and Narrow Algebraic Data}\label{sec:nbdata}
In dimension $1$ the Topological Data we use is pairs $(X,C)$ of  spaces $X$ with a set $C$ of base points. The Algebraic Data is groupoids, with the functor $\H$ sending
$(X,C)$ to the fundamental groupoid $$\pi_1(X,C)= \pi_1(X,X \cap C)$$ on the set $X \cap C$ of base points, as introduced in \cite{B67}. (This notation is a ruse to allow one to write $\pi_1(U,C)$ if, for example, $U$ is a subspace of $X$.) See Section \ref{sec:gpds}.

For modelling homotopy type in dimensions $> 1$  this method has to be refined owing to the
increased variety of convex sets in these dimensions. We then use at least two types of algebraic
data for a model;  their advantage is that while they are chosen to give equivalent categories,  these categories  serve
different purposes.  Thus in the following diagram, $\H$  of Criterion \ref{crit1}  has split into $\rho$ and $\Pi$.
Further $U$  is the forgetful functor, and $B = U \B$  is the ``classifying space"  of the Algebraic Data,
whose definition usually involves the Broad Algebraic Data:

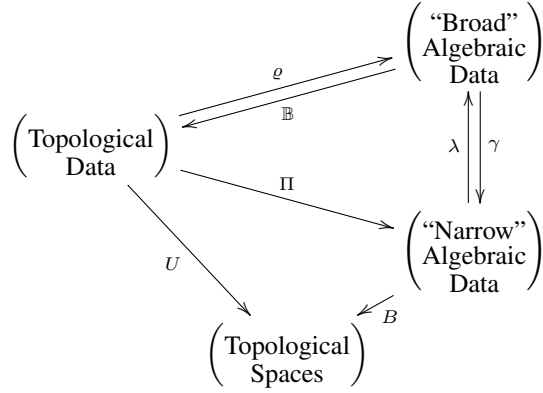
\begin{figure}[h]
\centering
 $\xymatrix @C=0.5pc@R=0.5pc{ &&\ar@<0.5ex>[dd]^{\gamma}   \left(\txt{``Broad''\\ Algebraic \\ Data}\right) \ar @<0.5ex>[dll]^{\mathbb B} \\
 \left(\txt{Topological \\ Data}\right)\ar @<0.5ex>[urr]^ \rho \ar[ddr]_U  \ar [drr]^\Pi  &&\\
  &&  \left(\txt{``Narrow'' \\Algebraic\\  Data}\right)\ar@<0.5ex>[uu]^{\lambda}  \ar [dl] ^ B \\ &\left(\txt{Topological \\ Spaces}\right)   &    }  $
 \caption{Broad and Narrow Algebraic Data} \label{fig:bddata}
\end{figure}
The ``Broad" Algebraic Data has a geometric type of axioms, is expressive,  and useful for conjecturing and
proving theorems, particularly colimit theorems, and for constructing classifying
spaces. Examples given later will be double groupoids, Section \ref{sec:hdim}, and cubical $\omega$-groupoids, Section \ref{sec:xmodclass}.

The ``Narrow" algebraic data has complicated algebraic axioms,  but relates to classical
theory, is useful for explicit calculation, and colimit examples have led to new algebraic constructions.
Examples given later will be crossed modules and crossed complexes.

The necessarily {\it algebraic}  proof of equivalence of these two structures allows one to get the best of both
worlds.  The hardest part of the equivalence is usually the natural equivalence $\lambda \gamma \simeq 1$, which shows that the broad data can be recovered from the narrow data. Also the functor $\gamma$ usually involves choices. The  whole relationship is a key to the power of the theory, and is important in developing
aspects of it. Intermediate models may also have specific uses.  Of course  the use of various geometric models such as discs, globes, simplices, cubes, for
defining invariants and constructions is standard in homotopy theory.

We will give more details in Section \ref{sec:fltsp} of the case when the Topological Data is that of
Filtered Spaces.  The Narrow Algebraic Data is then crossed complexes, which date back in
essence to Blakers, \cite{Bl48}, and were developed by Whitehead,  \cite{W49CHII}. The Broad Algebraic Data
is that of cubical $\omega$-groupoids, \cite{BH81algcub}.

Most other aspects are explained in  \cite{BHS},
to which this article gives an introduction.

Another example which we plan to explain  in Part II is the work in \cite{BL87,BL2},
where the Topological Data is that of $n$-cubes of pointed spaces, the Broad Algebraic Data
is that of cat-$n$-groups, later called cat$^n$-groups, and the Narrow Algebraic Data is that of crossed squares for $n=2$,  \cite[Section 5]{Lod82}, and in general crossed  $n$-cubes of groups, \cite{ESt}.
\begin{remark}\label{rem:induct}
  Many  categories  $\AD$ of Algebraic Data, either broad or narrow, have in our examples also a family of subcategories $\AD^n$ classified by dimension $n \geqslant 0$. and these come with forgetful functors $\Phi^n : \AD^n \to \AD ^{n-1}$. It usually turns out that $\Phi^n$ is a bifibration of categories, i.e. a fibration and an opfibration (or cofibration), and even has left and right adjoints. As spelled out in \cite[Appendix B]{BHS}, this helps the calculation of colimits in $\AD^n$ using their determination in $\AD^{n-1}$. For example, in order to calculate colimits of groupoids we build them from knowing colimits in the category of sets, i.e. the object sets of the groupoids. \qed
\end{remark}

The  Criteria of Section \ref{sec:phil} rule out on various grounds a number of other types of algebraic data
and models, such as simplicial groups, weak $\infty$-groupoids. See also the Introduction to \cite{Ellis-3-type}.

 \section{Six Anomalies in Algebraic Topology}\label{sec:anomal}
 It is useful to consider a subject in terms of {\bf anomalies}, i.e. aspects which do not seem to fit with the general tenor of the subject.

\begin{anomaly}\label{thm:anom1} The fundamental group is nonabelian, but homology and higher homotopy groups are abelian\footnote{The early 20th century workers in algebraic topology were hoping for higher dimensional nonabelian versions of the fundamental group. For this reason,  \v Cech's lecture  on higher homotopy groups  to the 1932 ICM at Z\"urich was not welcomed;
it was even felt these groups  could not be different from the homology groups. So only a small paragraph on this work appeared in the Proceedings of the ICM. See  \cite{jam}, and the Obituary Notice \cite{Alex}. }. \qed
\end{anomaly}
\begin{anomaly}\label{thm:anom2}
 The traditional  Seifert-van Kampen Theorem is a theorem about groups, but does not compute the fundamental group of the circle, which is {\bf  the}  basic example in algebraic topology. \qed
\end{anomaly}

\begin{anomaly}\label{thm:anom3} Traditional algebraic topology is fine with composing paths  but does not allow for the {\it  algebraic expression} of
$$\vcenter{\xymatrix@M=0pt@=1pc{\ar @{-} [rrrrr] \ar @{-} [dddd]&&&&&\ar @{-} [dddd]\\&&&&&\\
&&&&&\\
&&&&&\\
\ar @{-} [rrrrr] &&&&&}}\qquad \leftrightarrow \qquad   \vcenter{\xymatrix@M=0pt@=1pc{\ar
@{-} [rrrrr] \ar @{-} [dddd]&\ar @{-} [dddd]&\ar @{-} [dddd]&\ar
@{-} [dddd]&\ar @{-} [dddd]&\ar @{-} [dddd]
\\\ar @{-} [rrrrr]&&&&&\\
\ar @{-} [rrrrr]&&&&&\\
\ar @{-} [rrrrr]&&&&&\\
\ar @{-} [rrrrr]&&&&&}}$$
 From left to right gives {\it subdivision}.
  From right to left should give {\it composition}.  What we need for higher dimensional, nonabelian,  local-to-global problems is:

\hspace{10em}   {\it Algebraic Inverses to Subdivision.}   \qed

\end{anomaly}

\begin{anomaly}\label{thm:anom4}
In traditional homology theory,  for the  Klein Bottle diagram
    $$\xybiglabels \xymatrix@M=0pt@=3pc{   \ar @{} [dr]|\sigma &\ar @{-} [l]_a|(0.6)\tip \ar @{-} [d]|(0.6)\tip ^b\\
 \ar@{-}  [u]|(0.6)\tip^b & \bullet  \ar@{-}  [l]|(0.6)\tip^a }$$
    we  have to write $\partial \sigma= 2b$, not
$$\partial (\sigma)= a+b -a+b .$$
One can get more refined by working in the cellular operator chains of the universal cover which gives one (Whitehead, \cite{W49CHII},  Fox \cite{CF63}):
\begin{equation}\label{equ:bdyoper}
  \partial \sigma= a^{b-a+b} +b^{-a+b}-a^b +b.
\end{equation}
But this is clearly  more complicated,  and,  for subtler reasons,  less precise,   than the earlier formula\footnote{Actually the  formula \eqref{equ:bdyoper} has advantages, as shown in \cite{CF63}; it deals with modules over the group ring of a fundamental group. Nonetheless, information is lost, as shown in \cite[Section 7.4.v]{BHS}, and mentioned in \cite{W49CHII}.}. \qed
\end{anomaly}
\begin{anomaly} \label{thm:anom5}The product $\Delta^n \times \Delta^1$ of a simplex with the $1$-simplex, i.e. a unit interval, is not a simplex, although it has a simplicial subdivision, and this leads to awkwardness in dealing with homotopies in simplicial theory. \qed
\end{anomaly}
\begin{anomaly}\label{thm:anom6}
      There is no easily defined {\it multiple  composition} of simplices. \qed
\end{anomaly}

All the  Anomalies \ref{thm:anom1}--\ref{thm:anom6}  can be resolved  by using  in some way groupoids
 and their cubical  developments. While group objects in groups are just abelian groups, group objects in groupoids are equivalent to J.H.C. Whitehead's
 crossed modules, see  Definition \ref{defxmod},  which he used to describe the structure of  $$\pi_2(X,A,c) \to \pi_1(A,c), $$
 a major  example of a nonabelian structure in higher homotopy theory. This gives a functor $$\Pi: (\text{pointed pairs} ) \to (\text{crossed modules}).$$
How to compute it? We could do with  a Seifert-van Kampen type theorem!

\section{The origins of algebraic topology}\label{sec:orig}
The early workers in algebraic topology, before the subject had that name, see \cite{Di},  wanted to define numerical invariants using cycles modulo boundaries but
were not too clear about what these were!  Then Poincar\'e introduced formal sums of oriented
simplices and so the possibility of the equation $\partial \partial =0$, so that every boundary was a cycle. For more on the change from oriented to ordered  homology, see \cite{Bar95}.

The idea of formal sums of domains  came from an area with which many were concerned, see articles in \cite{jam},
namely integration theory, where one could write ``formally"

$$ \int _C f + \int _D f = \int_{C+D} f .$$
This link to integration theory automatically gives an abelian theory.

There was also a strong link to group theory through the study of Riemann surfaces.

In our account we use homotopically defined functors which involve actual compositions for the algebraic operations.  The use of
these in sufficient generality took much later developments in category theory, which itself led to familiarity with partial
algebraic operations. The aim was to find and use algebraic
structures which better model the geometry, and in particular the interactions of spaces.

\section{Enter groupoids}\label{sec:gpds}
I was led into this area in the 1960s through writing a text on topology\footnote{This was published as \cite{Brown-elements}, and is currently available as \cite{B2006}.}, and in particular on
the fundamental group $\pi_1 (X, c)$ of a space $X$ with base point $c$.  The standard Seifert-van Kampen Theorem  determines the fundamental group of a union of connected,  pointed
spaces. For example, the square of groups and morphisms induced by inclusions
\begin{equation} \vcenter{
\xymatrix{\pi_1(U \cap V,c) \ar [d] \ar [r] & \pi_1(V,c) \ar [d] \\
\pi_1(U,c) \ar [r] & \pi_1(X,c)}}
\end{equation}
is a pushout of groups if $U, V$   are open and $ c \in  U \cap V$, which is path connected.  A proof
by verifying the universal property, and more generally for unions of many sets, was given by
Crowell in \cite{Cr59}, following lectures by Fox. Note that the proof by the universal property  requires neither a previous proof that the category
of groups admits pushouts or colimits, nor  a specific construction of them.

If $U \cap V$   is not connected, there is no reasonable choice of base point, and of course a
standard example of this is the circle $S^ 1$.  A sensible resolution is to hedge your bets, not be
tethered to a single point, but use a set of base points appropriate to the geometry. Examples
of such choices are given in the following pictures, where the second depicts a union of three sets:

\begin{figure}[h]\centering
\includegraphics[width=4cm,height=2.5cm]{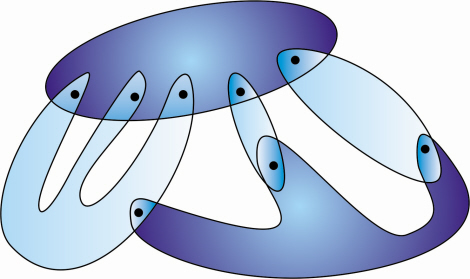}  \qquad \qquad  \includegraphics[width=3cm,height=3cm]{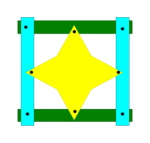}

\caption{Two unions of non connected spaces}\label{fig:unions}

\end{figure}
The right hand figure shows where a choice of one base point would destroy any symmetry
in the situation; this problem  is forcefully put in the paper on applications in partial differential equations,  \cite{PR15}. The left hand figure  is
used in \cite{Acb15}  to illustrate a failure for a space of both of two equivalent forms of the Phragm\'en-Brouwer Property\footnote{The reference \cite{Wilder1} proves the equivalence of these forms for connected and locally connected spaces.   One of the forms  is: if $D, E$ are disjoint closed subsets of $X$ such that neither $D$ nor $E$ separate $X$, then $D \cup E$ does not
separate $X$.}.

The following  Theorem was developed in order to resolve Anomaly \ref{thm:anom2}; it was published in \cite{B67}, and developed in the 1968, 1988, differently
titled, editions of \cite{B2006}.  It has also been presented in \cite{H71,Zi}.  Using it to relate  the Phragm\'en-Brouwer Property  to properties of fundamental groups  is a main step towards the proof  of the Jordan Curve Theorem in  \cite{B2006,Acb15}.
\begin{theorem} \label{thm:vktgpd} Let the space $X$ be the union of open sets $U,V$ with intersection $W$.
 The fundamental groupoid $\pi_1(X,C)$ on a set $C$ of base points, with the morphisms induced by inclusions,   gives a pushout of groupoids
 $$\xymatrix{\pi_1(W,C) \ar[d] \ar [r] & \pi_1(V,C) \ar [d]\\
 \pi_1(U,C) \ar [r] & \pi_1(X,C)}$$
 if $C$ meets each path component of $U,V,W$. \qed
\end{theorem}
More generally, the paper \cite{BRa84} considers a space $X$ with a cover by open sets $U_\lambda, \lambda \in \Lambda$ and a set $C$  of
``base points" which meets each path component of each $1, 2, 3$ -fold intersection of sets of the
cover. Then $ \pi_1(X,C)$  is given by a coequaliser of groupoids of the form
\begin{equation}\label{equ:equalisaer} \bigsqcup _{(\mu,\nu) \in \Lambda^2} \pi_1(U_\mu \cap U_\nu,C) \rightrightarrows \bigsqcup_{\lambda \in \Lambda} \pi_1(U_\lambda,C) \xrightarrow{c} \pi_1(X,C).   \end{equation}
Notice that an advantage of groupoids is that the coproduct of groupoids is simply disjoint union, so that the diagram \eqref{equ:equalisaer} is easily interpreted.
 The proofs of these theorems are direct generalisation of that in \cite{Cr59}, i.e. by verifying
the universal property, and so do not rely on knowing that pushouts or, more generally, colimits,  exist for groupoids, nor how to compute them. However the proof suggests that  groupoids are in some sense optimal for this theorem.

To model gluing you need to model spaces and maps.  But some of the maps may involve
identifications in dimension $0$, i.e. of points. This can be handled in the algebra of groupoids,
see \cite{H71} and \cite[Section 9.1]{B2006} and Remark \ref{rem:induct}, but clearly cannot even be formulated in an account which allows only path
connected spaces with one base point.

A basic example for such identification of points is to consider the map of pairs of spaces
$$ p: ([0,1],\{0,1\}) \to (S^1, \{1\}), $$
which just identifies $0$  and $1$. A consequence of Theorem \ref{thm:vktgpd} is that on taking $\pi_1$ we get the ``universal morphism'' of groupoids
$$\mathcal I\to \mathbb Z,$$ where $\mathcal I$  is the groupoid with two objects $0,1$  and exactly one arrow
$\iota : 0 \to 1$. The groupoid $\mathcal I$ is a generator for the category of groupoids, as is the infinite cyclic group $\mathbb Z $ for the category of
groups.

To see how this fits with the overall philosophy of Section \ref{sec:phil}, we take the Topological Data
to be a space $X$  and a set $C$; the Algebraic Data is groupoids; the functor $\mathbb H$  gives $\pi_1(X, C)$;
the functor $\B $ assigns to a groupoid $G$  the classifying space $B(G)$  together with the set $Ob(G)$.
Thus Anomaly \ref{thm:anom2} is resolved by the use of the fundamental groupoid on a set of base points.

See also \cite{B87} for wider uses of groupoids.

Theorem \ref{thm:vktgpd} gives a  transition from topology to algebra. The surprise to me of this kind
of result is that any particular group $\pi_1(X,c)$  is completely determined by the Theorem, but as
part of a much bigger structure which is in essence a model of the homotopy 1-type. So we have a process:
\begin{equation}
(X,C) = (\text{union}) \xrightarrow{\text{SvK Theorem}} \pi_1(X,C) \xrightarrow{(\text{combinatorics})} \pi_1(X,c) .
\end{equation}
However it may be that we do not want to reduce to some $\pi_1(X,c)$: for example a choice of a
$c \in C$  may obscure any symmetry of the situation.

One reason for the surprise is that exact sequences are a firm part of algebraic topology:
for adjacent dimensions they usually give information only up to extension, or even less.
That was the situation in the work with a Mayer-Vietoris sequence in nonabelian cohomology in \cite{B65}. Thus in this situation, an algebraic representative of the cohomology is more useful than the cohomology itself.

A possible reason for the success of such a representative seemed to be that groupoids have structure in dimensions
both $0$ and $1$, and this is surely needed to model homotopy $1$-types.  Philip Higgins told me of a dictum of
Philip Hall: ``You should seek an algebra which models the geometry and not try to force the geometry into a particular algebraic mode simply because it is more familiar."
Eldon Dyer had a dictum
that you should carry around the maximum structure as long as possible.  In particular,
throwing away the groupoid structure to get to a group is likely to be less productive. He also
suggested that I should take a hard line in pointing out that a connected space might be the
union of two open sets each with 50 path components and the number of path components of
the intersection being, say, 546.

Grothendieck wrote in \cite[24 April, 1983]{Gr-cor}:  \begin{quote}``....  the choice of a base point, and the $0$-connectedness assumption,
however innocuous they may seem at first sight, seem to me of a
very essential nature.  To make an analogy, it would be just impossible to work at ease with
algebraic varieties, say, if sticking from the outset (as had been customary for a long time)
to varieties which are supposed to be connected.   Fixing one point,  in this respect (which
wouldn't have occurred in the context of algebraic geometry) looks still worse, as far as limiting
elbow-freedom goes! Also, expressing a pointed 0-connected homotopy type in terms of a group
object mimicking the loop space (which isn't a group object strictly speaking), or conversely,
interpreting the group object in terms of a pointed ``classifying space", is a very inspiring magic
indeed  what makes it so inspiring it that it relates objects which are definitely of a very different
nature lets say, ``spaces" and ``spaces with group law". The magic shouldn't make us forget though
in the end that the objects thus related are of different nature, and cannot be confused without
causing serious trouble."
\end{quote}

All this suggests that for higher dimensions we should look for groupoid like algebraic
systems with structure in dimensions $0, 1, \ldots, n$.  This was one of the exciting possibilities I
considered in 1965, when it seemed that aspects of the proof of Theorem \ref{thm:vktgpd} could generalise
to the next dimension, see Section \ref{sec:hdim}.

That Section also explains how the use of groupoids suggested the possibility of resolving
Anomaly \ref{thm:anom1}.

\section{Higher dimensions?} \label{sec:hdim}
The abelian nature of higher homotopy groups to which we
referred in Anomaly \ref{thm:anom1} can be explained as follows:  the homotopy
classes of maps $(I^n , \partial I^ n ) \to (X, x)$  allows for $n > 1$ the structure of a set with more than $1$
compatible group structures, and a simple argument then shows that any two of these group
structures coincide and are abelian.

However it was known to the Grothendieck school in the 1960s, and was published independently
in \cite{BS1}, that a set with the two compatible structures of group and groupoid is not
necessarily abelian, and is equivalent to the structure known as a crossed module, as defined in
Definition \ref{defxmod}. We see later the occurrence of these structures in homotopy theory  through   second
relative homotopy groups; they also occurred in homological algebra, for representing elements
of $H^ 3 (G, M)$  where $G$  is a group and $M$  is a $G$-module.
\newpage
So the question could be formulated: are there higher homotopical invariants with structure
in dimensions from $0$ to $n$\footnote{The research method adopted was to look for obstructions to this prospect. If they were found, that could
be interesting. If they were, one by one, overcome, that would be even more interesting!\label{footn:meth}}? Are there colimit theorems in higher homotopy theory?
If so, what would they calculate?

The proof of the groupoid theorem seemed to generalise to dimension 2, at least, if one
had the right algebra of double groupoids, and the right gadget, a strict homotopy double
groupoid of a space. So this was an ``idea for a proof in search of a theorem".
Fortunately I found the book \cite{Ehresmann-65}: it contained the notion of double category and groupoid, which seemed the right concepts for resolving Anomaly \ref{thm:anom3}.

First we describe briefly the basic algebra of double groupoids, due to Charles Ehresmann.

We  think of a square in a space $X$ as a map $I^2 \to X$. Such a square is thought of as having four ``edges", each a path in $X$, as shown in the following picture.

\begin{equation}
\xybiglabels
\vcenter{\xymatrix@M=0pt@=3pc{\ar [d]_{\partial^-_2 } \ar [r]^{\partial^-_1} & \ar [d]^{\partial ^+_2}\\
\ar [r]_{\partial^+_1} &  } }  \qquad \directs{2}{1}
\end{equation}

Analogously,  in a double groupoid $G$  we have objects called squares and each square $\alpha$ has ``edges" $\partial^\pm_i$ for $i=1,2$, which lie in groupoids $G_1,G_2$, both with the  same set of objects $G_0$. One makes some obvious geometric conditions, so that a square has just 4 vertices.

Compositions in a double groupoid are allowed as follows: squares $\alpha, \beta$ have a composite $\alpha +_1 \beta$ if and only if $\partial^+_1 \alpha = \partial^-_1 \beta$, and the composition $+_1$ gives a groupoid structure.  Similarly we have a composition $+_2$.

\begin{alignat*}{3}
\xybiglabels \vcenter{\xymatrix@M=0pt{\ar@{-} [dd]\ar@{} [dr]|\alpha \ar@{-} [r] & \ar@{-} [dd] \\
\ar@{-} [r]\ar@{} [dr]|\beta& \\
\ar@{-} [r]& }}
&=
\begin{bmatrix}
  \alpha \\\beta
\end{bmatrix}& = \alpha+_1 \beta \\
\xybiglabels \vcenter{\xymatrix@M=0pt{\ar@{-} [d]\ar@{} [dr]|\alpha \ar@{-} [rr] &\ar@{} [dr]|\gamma \ar@{-} [d]& \ar@{-} [d]  \\
\ar@{-} [rr]&&  }} &= \begin{bmatrix}
  \alpha& \gamma
\end{bmatrix}& = \alpha+_2 \gamma \\
\xybiglabels \vcenter{\xymatrix@M=0pt{\ar@{-} [dd]\ar@{} [dr]|\alpha \ar@{-} [rr] & \ar@{} [dr]|\gamma\ar@{-} [dd]& \ar@{-}[dd]  \\
\ar@{-} [rr]\ar@{} [dr]|\beta&\ar@{} [dr]|\delta & \\
\ar@{-} [rr]&& }}  &= \begin{bmatrix}
\alpha & \gamma\\ \beta & \delta
\end{bmatrix} &\directs{2}{1}
\end{alignat*}
That each is to be a morphism for the other gives the interchange law:
 \begin{equation}\label{equ:interchange}  (\alpha\circ_2 \gamma)\circ_1(\beta\circ_2\delta)=
(\alpha\circ_1\beta)\circ_2(\gamma\circ_1\delta)
\end{equation}
\noindent for which the necessary and sufficient condition is that all the compositions in the equation are defined.  This equation illustrates  that a $2$-dimensional picture can be more comprehensible than a $1$-dimensional equation.

In these double groupoids the horizonal edges and vertical edges may come from different groupoids.

Note also that multiple compositions can be constructed  by considering ``compatible" arrays $(a_{ij}) $ of elements of a double groupoid, and their composite $[a_{ij}]$. Compatible means that each
$a_{ij}$ is composable with its neighbours.
 This easy expression of ``algebraic
inverses to subdivision" is one of the chief advantages of the cubical approach. It allows for the
potential of resolving Anomaly \ref{thm:anom3}, and also Anomaly \ref{thm:anom1}, since we will see later that double
groupoids can be much more complicated than abelian groups.  Indeed, classification of the
most general double groupoids is not well understood.
\section{Need also ``Commutative cubes"}\label{sec:commcub}
In dimension 1, we still need the 2-dimensional notion of commutative square:

\begin{equation}
\xybiglabels \vcenter{\xymatrix@M=0pt{\ar @{-}  [r]|\tip ^c\ar @{-}  [d]|\tip _a  &  \ar @{-}  [d]|\tip ^d \\
\ar @{-}  [r]|\tip _b & }} \qquad ab = cd ,   a = cdb^{-1}
\end{equation}

An easy result is that any composition of commutative squares is commutative.  In ordinary
equations:
$$ab = cd, ef = bg  \implies aef = abg = cdg.$$
The commutative squares in a groupoid G form a double groupoid $\square G$.

This fact is a key part of one proof of Theorem \ref{thm:vktgpd}, to verify that a morphism to a groupoid
is well defined. The essence is as follows: consider a diagram in a groupoid
\begin{equation}
\xybiglabels \vcenter{\xymatrix@M=0pt@=1pc{\ar @{-} [rrrrr]^b \ar @{=} [dddd]_1 &\ar @{-} [dddd]&\ar @{-} [dddd]&\ar @{-} [dddd]&\ar @{-} [dddd]& \ar @{=} [dddd] ^1 \\
\ar @{-} [rrrrr] &&&&&\\
\ar @{-} [rrrrr] &&&&&\\
\ar @{-} [rrrrr] &&&&& \\
\ar @{-} [rrrrr]_a &&&&& }}
\end{equation}
Suppose each individual square is commutative, and the two vertical outside edges are identities.
Then we easily deduce $a = b$.

For the next dimension we therefore expect to need to know what is a commutative cube,
and this is expected to be in a double groupoid in which horizontal and vertical edges come
from the same groupoid:  $G_1=G_2$.

\begin{equation}{{\objectmargin{0.1pc} \diagram &  \bu\rrto \xline '[1,0]
[2,0]|>\tip &&\bu\ddto\\\bu \urto \rrto \ddto &&\bu \urto \ddto &
\\&\bu \xline'[0,1] [0,2]|>\tip && \bu \\ \bu\rrto \urto &&\bu \urto
& \enddiagram }}\end{equation}

 We want the ``composed faces" to commute! What can this mean?

\begin{figure}[h]

$$  \xymatrix@M=0pt { & \ar@{..} [d]|a  \ar @{-} [r]  & \ar @{..} [d] |b & \\
 \ar @{-} [d] \ar @{..} [r] |a & \ar @{-}[r] \ar @{-}[d] & \ar @{-} [d] & \ar @{-} [d] \ar @{..} [l] |(0.25) {b^{-1}} \\
 \ar @{..} [r]|{c^{-1}} & \ar @{-} [r] \ar @{..} [d] |c & \ar @{..} [d] |{d}&  \ar @{..} [l] |{d} &\\
   & \ar @{-} [r] && } \qquad \xymatrix@M=0pt {\ar @{}[dr] |(0.35){\tl} & \ar@{..} [d]|a  \ar @{-} [r]  & \ar @{..} [d] |b &\ar @{}[dl]|(0.35)\tr \\
         \ar @{-} [d] \ar @{..} [r] |a & \ar @{-}[r] \ar @{-}[d] & \ar @{-} [d] & \ar @{-} [d] \ar @{..} [l] |(0.25) {b^{-1}} \\
         \ar @{..} [r]|{c^{-1}} & \ar @{-} [r] \ar @{..} [d] |c & \ar @{..} [d] |{d}&  \ar @{..} [l] |{d} &\\
         \ar @{}[ur]|(0.35)\bl & \ar @{-} [r] &&  \ar @{}[ul]|(0.35)\br }$$
\caption{``Composing" five faces of a cube}\label{fig:flatcomcub}
\end{figure}
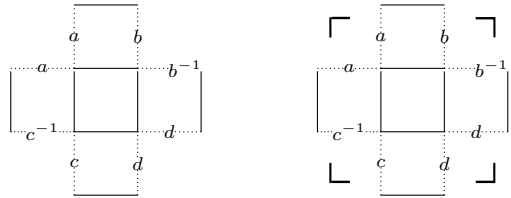

 We might say the top face is the composite of the other
faces:  so fold them flat to give the left hand diagram of Fig. \ref{fig:flatcomcub},   where the dotted lines show adjacent edges of a ``cut". We indicate how to glue these edges back together by extra squares which are a different kind of ``degeneracy" in the right hand diagram of this  Figure.

 Thus if we write the standard double groupoid identities in dimension 2 by
 $$ \tsq \quad \hh\quad \vv $$
 where a solid line indicates a constant edge, then the new types of square with commutative boundaries are written
$$ \tl \quad \tr \quad \bl \quad \br . $$
These new kinds of ``degeneracies" are called  {\bf connections}, \cite{BS76}\footnote{The latex code which was developed for this notation for  the identities and conections and used for example in \cite{BHS} is available at the following \url{http://www.groupoids.org.uk/connections-texcode.html} , Its advantages can be seen for interpreting diagram \eqref{equ:rotdiag}. The notation was first published in \cite{Brown-higherdimgroup}.}.
The laws on connections are:

\begin{equation} [\, \tl \;\br \;] = {\vv}\;\qquad
\begin{bmatrix}\,\tl\;\\ \br\,\end{bmatrix}  =
\hh\; \tag{cancellation}
\end{equation}

\begin{equation}
\begin{bmatrix}\tl & \hh \;\\ \vv & \tl
\end{bmatrix} = \tl\;\qquad
\begin{bmatrix}\br & \vv \;\\ \hh & \br  \end{bmatrix} = \br
 \tag{transport}
\end{equation}
 These are equations on turning left or right, and so
 {are a part of ``2-dimensional algebra". }

 The terms {\bf transport law} and  {\bf connection}
came from laws on path connections in differential geometry, see \cite{BS76}.

Any well defined $2$-dimensional composition of identities and connections  is called {\bf thin}. A thin square is entirely determined by its boundary.

Full details are given in \cite[Section 6.7]{BHS}.  This proof uses ``$2$-dimensional rewriting", as does the proof for strict cubical $\omega$-groupoids in Chapter 13 of \cite{BHS}. The account for the case of strict cubical  $\omega$-categories in \cite{higgins-thin} uses ``$3$-dimensional rewriting".

 Double groupoids allow for: multiple compositions, $2$-dimensional formulae,  and $2$-dimensional  rewriting.
As an example, we   get a ``rotation"

$$\sigma (\alpha)= \begin{bmatrix}
 \;  \vv & \tl & \hh\;  \\
  \bl & \alpha & \tr \\
  \hh & \br & \vv
\end{bmatrix}$$
As a taster to ``2-dimensional rewriting" we note that making two interpretations of the following diagram

\begin{equation}\label{equ:rotdiag}
\objectwidth{0in} \objectmargin{0in} \def\labelstyle{\textstyle}
\left[ \diagram
\midsq{\vv} & \midsq{\tl} & \midsq{\hh} & \midsq{\hh} & \midsq{\hh} & \\
\rline \midsq{\vv} & \rline \dline \midsq{\vv} & \rline \dline
\midsq{\sq} & \rline \dline \midsq{\sq} & \rline \dline \midsq{\sq} & \\
\midsq{\bl} & \dline \midsq{\alpha} & \dline \midsq{\tr} & \dline \midsq{\sq} & \dline \midsq{\sq} & \\
\midsq{\hh} & \dline \midsq{\br} & \dline \midsq{\vv} & \dline \midsq{\sq} & \dline \midsq{\sq} & \\
\midsq{\sq} & \dline \midsq{\sq} & \dline \midsq{\vv} & \dline \midsq{\tl} & \dline\midsq{\hh} & \\
\midsq{\sq} & \dline \midsq{\sq} & \dline \midsq{\bl} & \dline \midsq{\beta} & \dline \midsq{\tr} & \\
\midsq{\sq} & \dline \midsq{\sq} & \dline \midsq{\sq} & \dline \midsq{\vv} & \dline \midsq{\vv} & \\
\rline \midsq{\hh} & \rline \midsq{\hh} & \rline \midsq{\hh} & \rline \midsq{\br} & \rline \midsq{\vv} & \\
& & & & & \enddiagram \right] \xymatrix{ &\\ &\\ &\\
\ar@{}[d]_{.} & \\& }
\end{equation}
shows that  $\sigma\begin{bmatrix}
  \alpha &\beta
\end{bmatrix}= \begin{bmatrix}
   \sigma \alpha\\\sigma \beta
\end{bmatrix}$. See \cite[Theorem 6.4.10]{BHS}. Interpreting this proof in a homotopy double groupoid shows in principle how to construct explicit homotopies.    Diagrammatic arguments using the rules for connections are a central feature of that book.

It turns out that connections are what is needed,   \cite{BS76}, to obtain that  double groupoids with connections are equivalent to
crossed modules; these two structures  give respectively broad and narrow algebraic data.

This result was found
by trying to construct examples of double groupoids. The first example, given by Ehresmann, \cite{Ehresmann-65},
was that of commutative squares in a group or groupoid $G$.  The next example, say in the
case $G$ is a group, was to look at squares which commuted up to an element of a subgroup $N$,
say:  to get compositions one quickly found $N$ had to be normal in $G$.

  Finally, we looked at
a morphism $\nu : N \to G$ of groups and quintuples $(n; a, b, c, d)$ where $ n \in  N, a, b, c, d\in G $ and
$\nu n = b^{-1}a^{-1}cd$.
Compositions are easy to define if $G$ operates on $N$ satisfying CM1) of Definition \ref{defxmod}, but to get
the interchange law we need exactly also CM2). The thin squares are of the form $(1; a, b, c, d)$.
For more details, see also \cite{BrownMosaTAC}.

These ideas led to the notion of ``higher dimensional group theory" as a theory to be developed in the
spirit of group theory, but using algebraic structures involving partial operations
whose domains were defined by geometric conditions.   This scheme fitted with the work of
Higgins in \cite{Higgins-algebrawithoperators}.

\section{Cubical sets in algebraic topology}\label{sec:cubsets}
Dan Kan's first paper, \cite{Kan55}, was cubical, relying clearly on geometry and intuition.  It was
then found that cubical groups, unlike simplicial groups, were not Kan complexes, and there
was also a problem on realisation of cartesian products.  The simplicial methods seemed entirely adequate,
and the cubical theory seemed unfixable; so cubical methods were generally
abandoned for the simplicial, although many workers found them intuitive and useful in certain contexts. The book \cite{Mas80} gives a cubical approach to singular homology.

As explained in the previous Section, the work with Chris Spencer on double groupoids in the early 1970s, \cite{BS76}, found it
necessary to introduce an extra and new kind of ``degeneracy" in cubical sets, using the monoid structures
$$ \max, \min : [0, 1]^2 \to [0,1].$$

We called these ``connection operators". The  axioms for these in all dimensions were published in \cite{BH81algcub}.
Andy Tonks proved in \cite{tonks-cubgps}  that cubical groups with connections are Kan complexes!
G. Maltsiniotis in \cite{maltsin-cubandconn} showed that up to homotopy, connections correct the realisation
problem. Patchkoria in \cite{Pat12}  gave a cubical approach to derived functors. Cubical sets have
analogously been used for work on motives, see \cite{Vez14},  and the references there\footnote{Vezzani has explained to me the use as follows:
The cubical theory was better suited than the simplicial
theory when dealing with (motives of) perfectoid spaces in characteristic $0$.  For example: degeneracy maps of
the simplicial complex $\Delta$  in algebraic geometry are defined by sending one coordinate $x_ i$  to the sum of two
coordinates $y _j + y_{j+1}$.  When one considers the perfectoid algebras obtained by taking all $p$-th roots of the
coordinates, such maps are no longer defined, as $y _j + y_{ j+1}$  doesn't have $p$-th roots in general.  The cubical
complex, on the contrary, is easily generalized to the perfectoid world.}.

Cubical sets have also been used in concurrency theory, see \cite{goubault-geometric}:  if many computers run databases, then this can be envisaged as a system with many-dimensional time. They are also used in Homotopy Type Theory, \cite{CCHM}. A recent work, \cite{RZ16}, relates cubical sets with connections to the simplicial method of necklaces. The paper \cite{Lucas} investigates kinds of inverses in cubical higher categories with connections.

Independently of these facts, we deal with cubical sets with connections and compositions,
where the compositions have properties analogous to those of paths in a space; if all the compositions are groupoid structures, we get a cubical $\omega$-groupoid, \cite{BH81algcub}, \cite{BHS}.

\section{Strict homotopy double groupoids?}
In 1965, inspired by a definition in \cite{Ehresmann-65} which seemed to give a concept appropriate to resolve Anomaly \ref{thm:anom3},   I started to seek a strict homotopy double groupoid of a space.
Group objects in groups are abelian groups, by the interchange law.  Chris Spencer and
I found out in the early 1970s that group objects in groupoids are more complicated, in fact
equivalent to Henry Whitehead's crossed modules, see Definition  \ref{defxmod}, as was known earlier to some, see \cite{BS1}. The important point was that double groupoids could be highly non trivial\footnote{The further question was: could the leap to these higher dimensional versions of groups  be analogous to the leap from real analysis to many dimensional analysis? cf, Footnote \ref{footn:meth}.}.

In the early 1970s Chris Spencer, Philip Higgins and I developed a lot of understanding of: \\
(i) relations between double groupoids and crossed modules; and\\
(ii) algebraic constructions on the latter, including induced crossed modules, and colimit calculations.

In June, 1974,  Phil and I, stimulated by the last week of his supported visit to Bangor, argued as follows:

\begin{strategic}   \label{thm:strategy}

\begin{enumerate}[1.]
\item[]
\item As explained earlier, Whitehead had a subtle theorem (1941--1949), \cite{W41,W49CHII}, that
$$\pi_2(A \cup \{ e^2_\lambda \},A,c) \to \pi_1(A,c) $$
is a free crossed module.  So this was one, indeed seemingly the only, example of a
universal property in 2-dimensional homotopy theory.
\item Thus Whitehead's theorem should be a consequence of our conjectured theorem, if that so far unformulated theorem was to be any good.
\item But Whitehead's theorem was about second {\it relative}  homotopy groups.
\item So we should look for a homotopy double groupoid in a relative situation, $(X, A, c)$.
\item  So we tried what seemed the simplest idea, i.e. to look at maps of the square $I^ 2 \to X$  which took the
edges of the square to $A$  and the vertices to $c$, and consider homotopy classes of these.
\item  Because of all the preliminary work with Chris and Phil, this worked like a dream! It was
submitted 1975, and published severely cut as \cite{BH78sec}.
\end{enumerate}\qed

\end{strategic}

\section{Crossed Modules}\label{sec:xmod}
\begin{Def}\label{defxmod}   A  \emph{crossed  module} consists  of  a  morphism  of  groupoids  $\mu   :  M  \to  P$   of
groupoids such that $\mu$  is the identity on objects, $M$  is totally disconnected, and there is also
given an action of $P$  on $M$   written $(p, m)  \mapsto  \; ^p m$,  where if $p  :  x  \to y$   and $ m \in  M(x)$,
then $^p m \in M(y)$ .   The usual rules for such an action should be satisfied as well as for all
$ m, n \in  M, p \in  P $\\
CM1) $\mu (^pm )= p(\mu m)p^{-1}$, and\\
CM2)\footnote{To my knowledge, this rule was first stated in a footnote on p.422 of \cite{W41}. } $mnm^{-1} = \; ^{\mu m} n$.  \\
See Diagram \ref{eq:CM2} for a double groupoid explanation of this rule. \\
Crossed modules and their morphisms form a category $\XMod$.\qed
\end{Def}

If $P$   is a group, the notion of free crossed $P$-module $C(R) \to  P$  on a function $ w : R \to  P$
can be expressed by the pushout of crossed modules:
\begin{equation}\label{eqn:freexmod}
\vcenter{\xymatrix{\ar [r] \ar [d]   (1 \to F(R)) & \ar [d] (1 \to P)\\
\ar [r] (F(R)\xrightarrow{\;\; 1} F(R)) & (C(R) \to P)  ,   }  }
\end{equation}
where $F (R) \to F (R)$  is the identity crossed module.  Whitehead's Theorem solves the ``Klein
Bottle Anomaly"  \ref{thm:anom4}: the element $\sigma$  is the generator of $\pi_2(K ^2 , K^ 1 , x))$  as a free crossed $\pi_1(K^1,x)$-module.
This crossed module gives ``nonabelian chains" in dimension $\leqslant  2$. Whitehead's proof
of his theorem, see \cite{B80},  seems unlikely to suggest generalisations.

\section{Groupoids in higher homotopy theory?} \label{sec:gpdshhom}
Consider second relative homotopy groups $\pi _2 (X, A, c)$ .
\[\vcenter{\xybiglabels
\xymatrix@M=0pt@=3pc{\ar @{-}[r]^A \ar@{}[dr]|X \ar @{=}[d]_c \ar
@{=} [d] & \ar @{=}[d]^c
  \\
\ar @{=}[r]_c& }}\qquad \xdirects{2}{1} \]
where double lines show constant paths.   Note that the definition involves {choices},  and is unsymmetrical with respect to
directions.   This is unaesthetic! Further, all compositions are on a line:

$$\xymatrix@M=0pt@=2pc {\ar@{-}[rrrrrrrrr]\ar @{=}
[d] &\ar @{=} [d]  &\ar @{=}[d]   &\ar
@{=} [d]  & \ar @{=} [d] \ar
@{=} [d] &\ar @{=} [d]
&\ar @{=} [d] \ar @{=} [d] &\ar @{=} [d]
\ar @{=} [d] &\ar @{=} [d]
&\ar @{=} [d] \\
\ar@{=}[rrrrrrrrr]&&&&&& &&&} $$
All this is necessary to obtain a group structure, and, by tradition, the structure  has to be a group!

In 1974 Philip Higgins and the author considered a symmetrical version of relative homotopy groups, namely  $\rho_2(X,A,C)$, the   homotopy classes
rel vertices of maps $[0,1]^2 \to X$ with edges mapped into $A$ and
vertices into $C$:

 \[\vcenter{\xybiglabels
\xymatrix @=3pc{ \ar @{}[dr]|X \ar @{-}[r]|A  C
 \ar @{-} [d]|A    & \ar @{-}[d] |A  C \\ \ar @{-}[r]|A C& C}}\qquad \xdirects{2}{1}  \]

$$ \xymatrix{\rho_2(X,A,C)\ar @<0.33ex> [r]\ar @<-0.33ex>
[r]\ar @<1ex> [r]\ar @<-1ex> [r] & \pi_1(A,C)\ar @<0.5ex> [r] \ar
@<-0.5ex> [r]& C }
$$

 This allowed for the childish idea of gluing   two squares if the right side of one is the same as the left side of the other. This gives a partial   algebraic operation defined under geometric condition.

  There is a horizontal composition $$\langle \!\langle
\alpha \rangle \!\rangle +_2\langle \!\langle \beta \rangle
\!\rangle= \langle\!\langle \alpha +_2 h +_2 \beta \rangle\!\rangle$$ of classes in $\rho_2(X,A,C)$, where thick lines  show
constant paths.

$$\xybiglabels \vcenter{\xymatrix@=3pc@M=0pt{ \ar @{-}[r] \ar@{}[dr]|X \ar @{-}[d]& \ar @{}[dr]|A
\ar @{=}[r] \ar @{-}[d]&\ar @{}[dr]|X \ar @{-}[r]\ar @{-}[d]&\ar @{-}[d] \\
\ar @{-}[r]_\alpha & \ar @{=}[r]_h&\ar @{-}[r]_\beta& \\
&& }}\qquad  \directs{2}{1}
$$

 Intuition: gluing squares exactly is a bit too rigid; varying edges in $X$ is too free; but  ``varying edges in $A$'' seems just about right.

\begin{figure}[h]
\centering
  \includegraphics[width=3.8in,height=1.4in]{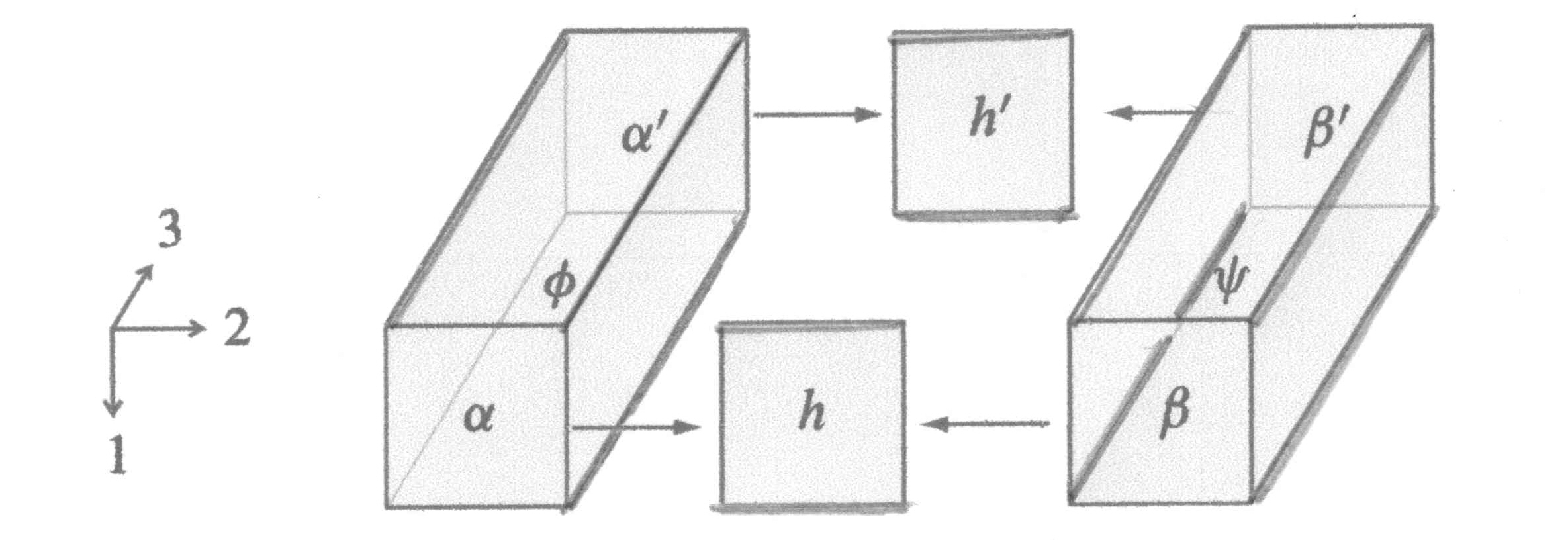}
\caption{Showing $+_2$ is well defined on $\rho$ }\label{fig:hole}
\end{figure}

 To show $+_2$ well defined,  let $\phi: \alpha \equiv
\alpha'$ and $\psi: \beta \equiv \beta'$, and let
$\alpha'+_2h'+_2 \beta'$ be defined. We get Figure \ref{fig:hole} in
which thick lines denote constant paths. Can you see why the `hole'
can be filled appropriately?

We find that $\rho(X, A, C)$ has in dimension 2 compositions in directions 1 and 2 satisfying
the interchange law and is a double groupoid with connections, containing as an equivalent
substructure the classical
$$\Pi(X.A,C) = (\pi_2(X.A,C) \to \pi_1(A,C)),$$
a crossed module over a groupoid. All this needs proof, but in this dimension is not too hard.

The second crossed module rule CM2) in Definition \ref{defxmod} is expressed in double groupoid terms as evaluating in two ways the following  diagram, where $\mu n =a$:
\begin{equation}\label{eq:CM2} \xybiglabels
  \vcenter{\xymatrix@M=0pt@=3pc{\ar @{-}[rrr] \ar @{=} [dd]\ar @{} [r] ^a \ar @{}[dr]|{\vv} &\ar @{=} [dd] \ar @{} [r]^{\mu m} \ar @{}[dr]|m &\ar @{=} [dd] \ar @{}[dr]|{\vv}\ar @{} [r] ^{a^{-1}} & \ar @{=}[dd]\\
  \ar @{=} [rrr] \ar @{}[dr]|{n} &\ar @{} [dr]|{\sq}& \ar @{}[dr]|{-_2n}& \\
  \ar @{=} [rrr]&& & }} \qquad \xdirects{2}{1}
\end{equation}

Thus $\rho$  gives in this case the ``broad" structure and $\Pi$  the ``narrow" structure.

Now we can directly generalise the 1-dimensional proof, since one has a homotopy double
groupoid which is ``equivalent"  to crossed modules but which can express:

\noindent $\bullet$ algebraic inverse to subdivision\\
$\bullet$ the notion of  commutative cubes, with the property  that any multiple composition of commutative cubes is commutative.

\begin{figure}[h]
\centering
\includegraphics[width=12cm,height=4cm]{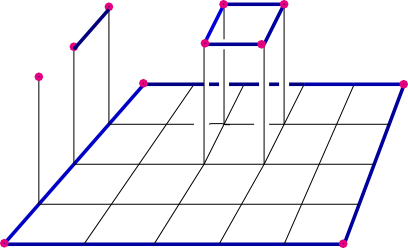}
\caption{Partial deformation of a subdivided square}\label{fig:deform}
\end{figure}

A key deformation idea is shown in Figure \ref{fig:deform}.
We need to deform the bottom subdivided square into a subdivided square for which all of its
subsquares define an element of $\rho(X, A, C)$. This explains the connectivity assumptions for the
Theorem. It is to express this diagram that $\rho(X, A, C)$ is designed.

We end up with a 2-d SvK theorem, namely a pushout of crossed modules:
\begin{equation}\label{equ:2push}
  \vcenter{\xymatrix{ \Pi(X,A,C) \cap U \cap V \ar [r] \ar [d] & \Pi(X,A,C) \cap V \ar [d] \\
  \Pi(X,A,C) \cap U \ar [r] & \Pi(X,A,C) }}
\end{equation}
if $X = U \cup V$, $U, V$  are  open, and the intersections of $(X,A,C)$ with $U,V, U \cap V$ are  connected in the following sense:

\begin{Def}\label{def:connpair} We say $(X, A, C)$  is  \emph{connected} if \\
(i) $\pi_0(C) \to \pi_0(A), \pi_0(C) \to \pi_o(X) $ are surjective; and \\
(ii) any map $(I, \partial I) \to (X, C)$ is deformable into $A$ rel end points. \qed
\end{Def}

The huge generalisation of Whitehead's Theorem given by the pushout of diagram \eqref{equ:2push}, says that you can glue homotopy 2-types, by
gluing crossed modules as a pushout.
 A consequence is the following pushout, for $A,X$  connected and $f: A \to X$:
\begin{equation}
\vcenter{\xymatrix{\ar [r] \ar [d]   (1 \to \pi_1(A,a)) & \ar [d] (1 \to \pi_1(X,f(a)))\\
\ar [r] (\pi_1(A,a) \xrightarrow{\;\; 1} \pi_1(A,a)) & (\pi_2(X \cup_f C(A),X,f(a))  \to \pi_1(X,f(a)) ) .   }  }
\end{equation}
Whitehead's Theorem is the case $A$ is a wedge of circles.

The results enable some nonabelian computations of homotopy 2-types, see \cite[Chapters 4,5]{BHS}.  There are methods to determine explicitly the module $\pi_2$  from this algebraic
model of the 2-type, but that module is a pale shadow of the 2-type.

The proof of all this works first in the category of double groupoids with connection, and then uses the equivalence with crossed modules given in \cite{BS76}.  The notion of connection is the precisely  necessary tool   for both parts of the proof.

To define the functor $\B$ on crossed modules it is best to generalise to crossed complexes and their relation with filtered spaces.
\section{Filtered spaces}\label{sec:fltsp}
The success of the 2-d work  led to a search for the $n$-dimensional idea, and in view of the work of \cite{Bl48} and \cite{W49CHII} it  was natural to look at filtered spaces:
\begin{equation}
  X_* := X_0 \subseteq X_1 \subseteq \cdots X_n \subseteq \cdots  \subseteq  X
\end{equation}
With the obvious morphisms, these form a category $\FTop$.  Examples of filtered spaces are:\\
(i) the skeletal filtration of a CW-complex: e.g. $\Delta^n_*, I^n_*$. \\
(ii) $\bullet \subseteq A \subseteq \cdots \subseteq A \subseteq X \subseteq X \subseteq \cdots $,
where the first $X$ occurs at position $n$;\\
(iii)  for a pointed topological space $X$,  the free monoid $F M(X)$ on $X$ may be filtered by word length, cf. \cite{Ba-B}.

Now we need to define the notion of a cubical set with compositions and connections and
then check that for a filtered space $X_*$  we get an example of this if we define $R(X _* )$  by
$$ R_n(X_*) = \FTop(I^n_*, X_*).$$
Of course the compositions here do not form groupoids.

Now define a {\it thin} homotopy $h: a \equiv b$ between elements $a,b \in R_n(X_*)$  to be a homotopy through filtered maps and rel vertices of the cube.  Such a homotopy is a particular kind of element of $R_{n+1}(X_*)$. The thin homotopy class of an element $a \in R_n(X_*)$ is written $[a]$.
\begin{theorem} \label{thm:rho} The thin homotopy classes of elements of $R(X_*)$ form a cubical set with connections written $\rho(X_*)$ with projection
$p: R(X_*)  \to \rho(X_*)$. Then: \\
(i) the compositions on $R(X_*)$ are inherited to give $\rho(X_*)$ the structure of strict cubical $\omega$-groupoid;\\
(ii) $\rho(X_*)$ is a Kan complex in which every box has a unique thin filler,  where an element
$\alpha  = [a] \in \rho_n(X_*)$ is  \emph{thin} if it has a representative $a$  such that $a(I^ n ) \subseteq X_{ n -1}$;\\
(iii) the projection $p$ is a Kan fibration of cubical sets.  \qed
\end{theorem}
The proof was given in \cite{BH81:col}, and full details are also in \cite{BHS};   the proof of (i) uses cubical versions of Whitehead's notions of expansions and collapsing,
and it seems significant  that there is exactly enough ``filtered room"  to make the proof work.

The proof  that $R(X_*)\to  \rho(X_*)$ is a fibration, relies on a nice
use of geometric cubical methods, and the Kan condition.  One applies the Kan condition by
modelling it in subcomplexes of real cubes and using expansions and collapsings of these.

A consequence of the  fibration property is the following, which we give here because
it emphasises the use of multiple compositions:

\begin{cor}[Lifting composable arrays]  Let $(\alpha_{(i)} )$  be a composable array of elements
of $\rho_n(X_*)$. Then there is a composable array $(a_{(i)})$  of elements of $R_ n(X_* )$  such that for all $(i)$,
$p(a_{ (i) } ) = \alpha _{(i) }$. \qed
\end{cor}
Thus the weak cubical infinity groupoid structure of $R(X _* )$  has some kind of control by the
strict infinity groupoid structure of $ \rho(X_*  )$.

\section{ Crossed complexes and their classifying space}\label{sec:xmodclass}
Associated functorially with a filtered space $X_*$  is a homotopical invariant $\Pi(X_* )$  called a {\it crossed
complex}, due in the case $X_0 $  is a point to Blakers, \cite{Bl48}, and Whitehead, \cite{W49CHII}. It involves
the crossed module of groupoids
$$\pi_ 2 (X_ 2 , X_ 1 , X_0) \to \pi_1 (X _1 , X_0), $$
and in dimensions $> 2$  the boundary map
$$\pi_n (X_ n , X_ {n-1} , x) \to \pi_{ n-1} (X_{ n-1} , X_{ n-2} , x), x \in X_0$$
together with the operation of $\pi_ 1 (X _1 , X_ 0 )$  on all these groups. So we have a functor $\Pi :\FTop  \to \Crs$.

For our purposes, a key result is that the algebraic equivalence of categories $$\gamma  : \omega\text{-}\Gpd \to  \Crs$$  takes
$\rho(X _* )$  to $\Pi(X _* )$. So we may use whatever is best in a given situation. It also confirms that the axioms for crossed complexes and
for $\omega$-groupoids are the right ones. Thus the algebra matches the geometry.

Now we can define the classifying space $B(C)$ of a crossed complex $C$. The cubical definition
fits better with other themes. So we first define  the cubical set $N(C)$ in dimension $n$  by
$$N(C)_n= \Crs(\Pi(I^n_*),C).$$
We then define
$$B(C) = |N(C)|,$$
the geometric realisation of the cubical set $N(C)$.

But $C$ can be filtered by its truncations $C^ {(n)}$, the crossed complex which  agrees with $C$  up to and including
dimension $n$ and above that is trivial. So we get a filtered space $ \mathbb B (C)$ which in dimension $n$ is
$B(C^{ (n)})$. The proof of connectivity of $\mathbb B (C)$ uses the adjointness
$$\Cub(K, N(C)) \cong \Crs(\Pi(|K|_* , C)).$$

We also need to discuss connectivity of a filtered space $X_*$. The motivation of the condition is to allow the deformations used  in the proof of the Seifert-van Kampen Theorem, namely that the functor $\rho$, and hence also $\Pi$, preserves certain colimits.
\begin{Def}\label{Def:connfiltered}
  A filtered space $X_*$ is {\it connected} if the following conditions hold: for all $n > 0$ \\
  ($\phi 0$): the induced function $\pi_0 X_0 \to \pi_0 X_n$ is surjective and, \\
  ($\phi n$): $\pi_j(X_{n+1},X_n,v) =0 $ for all $j \leqslant  n$, and $v \in X_ 0$. \qed
\end{Def}
There are other equivalent forms of this definition, obtained using  homotopy exact sequences. It is standard that these conditions
are satisfied if $X_*$ is the skeletal filtration of a CW-complex, but this is also a consequence of the following

\begin{theorem}\cite{BH81:col}\label{thm:hhvkt} Let $X_*$ be a filtered space and let $\mathcal U = \{U^\lambda, \lambda \in \Lambda \}$ be an open cover of $X$. For each $\nu \in \Lambda ^n$ let $U^\nu$ denote the intersection of the sets $U^{\nu_i}$ and let $U^\nu _*$ be the filtered space formed by intersection with the $X_i, i \geqslant 0$. Suppose that for all $n \geqslant 1$ each  $U^\nu_*, \nu \in \Lambda ^n$  is a connected filtered space.
Then\\
\emph{(Con)} $X_*$ is connected, and \\
\emph{(Iso)} the following diagram is a coequaliser diagram of crossed complexes:
$$\bigsqcup _{\nu \in \Lambda ^2} \Pi U^\nu _*  \rightrightarrows^a_b \bigsqcup _{\lambda \in \Lambda } \Pi U^\lambda _* \xrightarrow{c}  \Pi X_* .\qed  $$
  \end{theorem}
Here $a,b$ are induced by the inclusions $U_{(\lambda,\mu)} \to U_\lambda, U_{(\lambda,\mu)} \to U_\mu$ respectively, and $c$ is induced   by the inclusions $U_\lambda \to X$.

This Higher Homotopy Seifert-van Kampen Theorem allows one easily to compute for example $\Pi(X_*)$ when $X_*$ is the skeletal filtration of a CW-complex.

\section{Why filtered spaces?}\label{sec:whyfilt}
In the proof of   Theorem \ref{thm:rho}, particularly the proof that compositions are inherited, there
is exactly the right amount of ``filtered room".  One then needs to evaluate the significance of
that fact!

Grothendieck in ``Esquisse d'un Programme", \cite[Section 5]{GrEsq}, has attacked the dominance of the
notion of topological space, which he says comes from the needs of analysis rather than geometry.
He advocates some ideas of stratified spaces: filtered spaces are a step in that direction.
A general argument is that to describe or specify a space you need some kind of data, and
that data has some kind of structure. So it is reasonable for the invariants to be defined using
that structure.

The theory given so far gives an account of algebraic topology on the border between
homology and homotopy but without the use of singular homology or simplicial approximation.
Nonetheless it allows a geometric account of results such as:
\begin{enumerate}[(i)]
\item  the Brouwer degree theorem (the $n$-sphere $S ^n$  is $(n - 1)$ -connected and the homotopy
classes of maps of $S^ n$  to itself are classified by an integer, called the degree of the map);
\item   the Relative Hurewicz theorem, which usually relates relative homotopy and homology
groups, but is seen here as describing the morphism $$\pi_n(X, A, x) \to \pi_n(X \cup CA, CA,x) \xrightarrow{\cong} \pi_n (X \cup  CA, x)$$  when
$(X, A)$ is $(n - 1)$-connected;
\item   Whitehead's theorem (1949) referred to earlier;
\item   a generalisation of that theorem to describe the crossed module
$$ \pi_2(X \cup_ f CA, X, x) \to  \pi_ 1 (X, x)$$
as induced by the morphism $f_* : \pi_1(A, a) \to  \pi_ 1 (X, x)$  from the identity crossed module
$\pi_ 1 (A, a) \to \pi _1 (A, a)$; and
\item  a coproduct description of the crossed module $\pi_ 2 (K \cup L, M, x) \to  \pi_1 (M, x)$  when $M =
K \cap L$ is connected and $(K, M), (L, M)$ are $1$-connected and cofibred.
\end{enumerate}
The Fibration Theorem, part (iii) of Theorem \ref{thm:rho}, yields  an associated family of  homotopy exact
sequences which, in the case $X_*$  is a CW-filtration, corresponds to Whitehead's exact
sequence in \cite{W50}, as shown in \cite[Theorem 14.7.9]{BHS}.  Thus the base $\rho(X_*)$ of the fibration, captures
 the ``linear" part of the space $X$, while $R(X_* )$  has the homotopy type of $X$. Note also that the two papers  \cite{BH81algcub,BH81:col} take under 60 pages to cover items (i) - (iv) above.

The above cited book  gives applications of monoidal closed categories to homotopy classification
results, based on the key isomorphism $$I^ m_* \otimes I^n_* \cong I^{ m+n}_*. $$ For example it is easy to define in the
cubical setting a natural transformation
$$\eta  : \rho(X_* ) \otimes \rho(Y _* ) \to \rho(X_*\otimes  Y_* ),$$
basically as $[f] \otimes  [g] \mapsto  [f \otimes g]$, and then to transfer this to the crossed complex case using the
equivalence of categories. This leads to the notion of ``crossed differential algebra", i.e. monoid objects with respect to $\otimes$
in the category $\Crs$, though the theory may work out better in the category \ogpd.

A further advantage of the monoidal closed property of the category $\Crs$ is to yield  a  new view of the Homotopy Addition Lemma (HAL),  or Theorem: this gives an explicit formula for $\Pi(\Delta^n_*)$.  This result is used in \cite{Bl48}, was first proved in \cite{Hu53},  is proved in \cite{WG78} inductively  in conjunction with the Relative and Absolute Hurewicz Theorems. By contrast, in \cite[Theorem 9.9.4]{BHS} the HAL  is proved  in the category $\Crs$ by regarding the filtered simplex  $\Delta^n_*$ as a cone on $\Delta^{n-1}_*$ and using the detailed algebraic description of the tensor product of crossed complexes.

\begin{figure}[h]\centering
 \includegraphics[width=2in,height=2in]{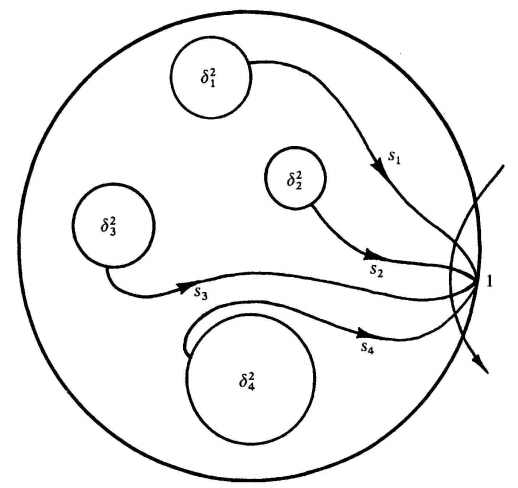}

\caption{2-dimensional Homotopy Addition Lemma } \label{fig:2dHAL}
\end{figure}

Interestingly, a kind of 2-dimensional form of the HAL plays a key role in Whitehead's proof of his theorem on free crossed modules: see Figure \ref{fig:2dHAL}, taken from  \cite[p.148]{B80},
which is used to describe an element of $\pi_2(X \cup \{ e^2_\lambda \},X,x)$ as given by a product of elements of the form $(c_{i_k}^{\xi_k})^{\zeta_k}$.

The crossed complex model of homotopy types  can be described as essentially a ``linear" model:  it does not even
capture quadratic information such as Whitehead products. Grothendieck was led by this limitation of modelling of homotopy types
to embark as one thread in ``Pursuing Stacks", \cite{GrPS},  on the notion of weak $\infty$-groupoids, on which
there is now a considerable body of  work\footnote{ See \url{https://ncatlab.org/nlab/show/infinity-groupoid}. This literature currently makes little reference to  cubical methods. }.

The linearity limitation can also be seen intuitively in the fact that crossed complexes can be  written on a line, and crossed resolutions are constructed step by step in increasing dimensions.

By contrast, Loday's construction of a ``resolution" of a connected, pointed space in \cite{Lod82},
corrected in \cite{St}, moves a step in one direction, then in an orthogonal direction, and so on, to
end up with the Topological Data of connected $n$-cube of fibrations, and with Broad Algebraic
Data that of cat$^n$-groups, which are essentially $(n + 1)$-fold groupoids in which one direction,
say the last, is a group. The Narrow Algebraic Data is the crossed $n$-cubes of groups of \cite{ESt}. More on this will be written in a second paper.

One aspect from this work which has been well taken up is the notion of nonabelian tensor
product of groups\footnote{I have kept up a bibliography on this area, available as \url{www.groupoids.org.uk/nonabtens.html},
with currently 156 items dating from 1952. Most interest has been from group theorists. }, defined in \cite{BL87}.  Nonetheless, this restriction to pointed spaces is
yet another Anomaly!  Further, it can be expected that quadratic information can be recovered by using crossed bicomplexes, i.e. crossed complexes internal to the category $\Crs$.

The Algebraic Model of $n$-fold groupoids is considered in \cite{BB12}, but currently this work has not   been linked with the work on higher order van Kampen theorems, i.e.  following the philosophy of Section \ref{sec:phil}.

Contact: \url{www.groupoids.org.uk}

\end{document}